\documentclass[]{article}
\usepackage{amsthm}
\usepackage{amsmath}
\usepackage{amssymb}
\usepackage{graphicx}
\usepackage{epstopdf}

\newtheorem{theorem}{Theorem}
\newtheorem{lemma}{Lemma}

\newtheorem*{theorem*}{Theorem}
\newtheorem*{claim}{Claim}

\theoremstyle{definition}

\newtheorem*{nota}{Notation}
\newtheorem{step}{Step}

\theoremstyle{remark}
\newtheorem*{rem}{Remark}

\newtheorem{cond}{Condition}

\title{Self-similar solutions for dyadic models of the Euler equations}
\author{In-Jee Jeong}
\date{\today}

\global\long\def\g{\gamma}
\global\long\def\G{\Gamma}
\global\long\def\D{\Delta}

\global\long\def\l{\lambda}
\global\long\def\V{\Vert}
\global\long\def\pr{\partial}

\usepackage{geometry}
\begin{document}

\maketitle

\begin{abstract}
	We show existence of self-similar solutions satisfying Kolmogorov's scaling for generalized dyadic models of the Euler equations, extending a result of Barbato, Flandoli, and Morandin \cite{MR2746670}. The proof is based on the analysis of certain dynamical systems on the plane.
\end{abstract}

\section{Introduction}

In this note, we address the question of the existence of \textit{self-similar} solutions in the following infinite system of ordinary
differential equations:
\begin{equation}\label{eq:KPO}
\begin{split}
\frac{da_{j}(t)}{dt} &= \big(\lambda^{j}a_{j-1}^{2}(t)-\l^{j+1}a_{j}(t)a_{j+1}(t)\big)\\
& \qquad +\beta\big(\l^{j}a_{j-1}(t)a_{j}(t)-\l^{j+1}a_{j+1}^{2}(t)\big),
\end{split}
\end{equation}
for $j\geq0$, with the boundary condition $a_{-1}(t)\equiv0$.

The factor $\l$ is some fixed constant greater than 1, and the coefficient $\beta$ is taken to be a nonnegative constant. The case 
$\beta=0$ is sometimes called the KP equations, and have first appeared
in the literature in the work of Friedlander, Katz, and Pavlovic \cite{MR2038114,MR2095627}. The full system \eqref{eq:KPO} was suggested in a work of Kiselev and Zlatos \cite{MR2180809} -- it was characterized as an infinite system of ODEs which is quadratic, conserves the energy, and contains only the nearest neighbor interactions (see \cite[Proposition 2.4]{MR2180809}). It was shown in \cite{MR3339169} that smooth solutions of \eqref{eq:KPO} blows up in finite time when $\beta $ is small enough, extending previous results which established blow-up in the case $\beta = 0$ \cite{MR2038114},\cite{MR2095627},\cite{MR2180809},\cite{MR2231615},\cite{MR2415066}. 
This type of equations are suggested as toy models for the dynamics of an inviscid fluid. Roughly speaking, the square of the scalar variable $a_j(t)$
is associated with the energy of a fluid velocity vector field restricted to a frequency shell of radii $\sim 2^j$. The quantity $\sum_{j \ge 0} a_j^2(t)$
is then the analogue of energy, and one may check directly from \eqref{eq:KPO} that it is (formally) conserved in time.\footnote{It is justified, for example, when $\sum_{j\ge 0} \l^{2j/3+\delta}a_j^2(t) < \infty$ for some $\delta>0$.} Therefore, it is quite natural to restrict to non-negative and finite-energy solutions, i.e. a sequence of functions $a_j(t) \ge 0$ solving \eqref{eq:KPO} with $\sum_{j \ge 0} a_j^2(t)<\infty$. 

If one attempts to consider \eqref{eq:KPO} as a model for \textit{turbulence}, it makes sense to add a constant forcing term to the lowest mode (to sustain turbulent motion): 
\begin{equation}\label{eq:KPO2}
\begin{split}
\frac{da_{j}(t)}{dt} &= \big(\lambda^{j}a_{j-1}^{2}(t)-\l^{j+1}a_{j}(t)a_{j+1}(t)\big)\\
& \qquad +\beta\big(\l^{j}a_{j-1}(t)a_{j}(t)-\l^{j+1}a_{j+1}^{2}(t)\big) + f \delta_0(j) ,
\end{split}
\end{equation}
where $f> 0 $ is a constant and $\delta_0(j) = 1$ for $j = 0$ and 0 otherwise. In the case $\beta = 0$, it is elementary to check that there exists a unique (finite-energy) fixed point\footnote{The existence of a fixed point in a forced system contradicts energy conservation -- this phenomenon is called either \textit{anomalous} or \textit{turbulent} dissipation.} of the system \eqref{eq:KPO2}, and it has the form \begin{equation}
a_j(t) = \mathrm{const}(f)\cdot \l^{-j/3}.
\end{equation}
It can be argued that the scaling $\l^{-j/3}$ corresponds to Kolmogorov's famous 5/3 law (see \cite{MR2337019},\cite{MR2600714} for example). Remarkably, it was shown that this fixed point attracts all other solutions \cite{MR2600714} when $ t \rightarrow +\infty$.\footnote{Solutions exist globally in $l^2$, even though they blow up in finite time with respect to smooth topologies.} Let us note that a fixed point with the same scaling $\l^{-j/3}$ continues to exist for $\beta > 0$. 


In the absense of forcing\footnote{This may be viewed as a model for \textit{freely decaying turbulence}, see \cite{MR2746670}.}, nontrivial fixed points do not exist anymore, but there are special solutions which play a similar role (at least conjecturally). We say that a solution $a(t) = \{ a_j(t) \}_{j\ge 0}$ is \textit{self-similar} if there exists some profile $\phi(t)$ such that \begin{equation}\label{eq:SS}
a_j(t) = a_j^* \cdot \phi(t)  \quad\mbox{for all}\quad j \ge 0,
\end{equation}
with constants $a_j^*$. By plugging in \eqref{eq:SS} into \eqref{eq:KPO}, one readily sees that $\phi(t)$ must take the form $(t-t_0)^{-1}$ for some $t_0 \in \mathbb{R}$. In addition, when $\beta = 0$, the constants $a_j^*$ should satisfy the recurrence
\begin{equation}\label{eq:recursion_og}
a^*_{n+1} = \frac{1}{\l^{n+1}} +  \frac{(a^*_{n-1})^2}{\l a_{n}^*}.
\end{equation}
It is convenient to renormalize the variables by $\l^n a^*_n = \alpha_n$. Then \eqref{eq:recursion_og} takes the form
\begin{equation}\label{eq:recursion}
\alpha_{n+1} = 1 + \lambda^2 \frac{\alpha_{n-1}^2}{\alpha_{n}}.
\end{equation} In the case of $\beta > 0$, the corresponding recursion takes the form 
\begin{equation}\label{eq:recursion_plus}
\alpha_{n+1} = -\frac{\l}{2\beta} \alpha_n + \sqrt{ \big( \frac{\l}{2\beta} \alpha_n \big)^2 + \frac{\l}{\beta} \big( \l^2 \alpha_{n-1}^2 + \beta \l \alpha_{n-1} \alpha_n + \alpha_n \big)  }.
\end{equation}

In principle, any choice of $\alpha_0 > 0$ would yield a non-negative self-similar solution via \eqref{eq:recursion}, but we are interested in finite-energy solutions. It is surprising that this condition uniquely selects a value of $\alpha_0$:
\begin{theorem*}[Barbato-Frandoli-Morandin \cite{MR2746670}]\label{thm:BFM}
	There exists a unique value of $\alpha_0 > 0 $ such that the self-similar solution \begin{equation}\label{eq:selfsimilarsolution}
	a_j(t) = \frac{a_j^*}{t-t_0}
	\end{equation}
	obtained from the recursion \eqref{eq:recursion} and the scaling  $\l^j a^*_j = \alpha_j$ has finite energy. Moreover, we have asymptotics \begin{equation*}
	a_j^* \approx \mathrm{const} \cdot \lambda^{-j/3} \quad\mbox{as}\quad j \rightarrow \infty.
	\end{equation*}
\end{theorem*}

Our result extends the existence statement to the case of small $\beta > 0$. 

\begin{theorem}\label{thm:main}
	There exists $\beta_0 > 0$ such that for all $\beta \in [0,\beta_0)$, there exists a value $\alpha_0 = \alpha_0(\beta) > 0$ such that the sequence of points $\{\alpha_j^*\}_{j \ge 0}$ obtained from the recursion \eqref{eq:recursion_plus} and the scaling  $\l^j a^*_j = \alpha_j$ satisfies
	\begin{equation*}
	a^*_j \approx \mathrm{const}(\beta) \cdot \l^{-j/3} \quad\mbox{as}\quad j \rightarrow \infty.
	\end{equation*}
\end{theorem}

\begin{rem}
	The proof of \cite{MR2746670} was based on complex analysis and it is not clear to us whether the method can be adapted to the case $\beta > 0$. 
	
	Our arguments yield the uniqueness statement for $\beta =0$ as in \cite{MR2746670}, and also ``local'' uniqueness for $\beta > 0$ small (in the sense that if we slightly perturb $\alpha_0(\beta)$ a little bit, it does not provide a finite energy self-similar solution). 
\end{rem}

\begin{rem}
	As it was pointed out in \cite{MR2746670}, the unique solution given in the theorem automatically generates a family of self-similar solutions parametrized by $(J,t_0) \in \mathbb{Z}_+ \times \mathbb{R}_+$, where $t_0$ is the time parameter as in \eqref{eq:selfsimilarsolution} and $J$ is the first nonzero entry in the sequence. It is reported in \cite{MR2746670} that at least numerically, any solution of \eqref{eq:KPO} (with $\beta = 0$) \textit{selects} one of the self-similar solutions and converges to it exponentially in time. We have observed a similar phenomenon for the case of small positive $\beta$. 
\end{rem}

In the proof, we fix $\l =2$ but  the proof carries over to any value of $\l > 1$. Now let us present an outline of the proof of the Theorem \ref{thm:main}, for the case $\beta = 0$.

\begin{enumerate}
	\item The starting point is a reformulation of the statement in terms of a dynamical system defined on the plane. If we consider the map \begin{equation}
	F:(x,y) \mapsto (y, 1 + \l^2 \frac{x^2}{y} ),
	\end{equation}
	then the equivalent problem is to find a point $(0,\alpha_0)$ with $\alpha_0 > 0$ whose iterates under $F$ have the desired asymptotics. Note that $F$ restricts to a dynamical system in the positive quadrant $\mathbb{R}_+ \times \mathbb{R}_+ = \{ (x,y) : x, y >0 \}$ and is injective. As the image of the ray \begin{equation}
	L = \{ (x,y) : x = 0, y >0 \}
	\end{equation} is already contained in $\mathbb{R}_+ \times \mathbb{R}_+$, we may consider $F$ only on the positive quadrant. We then change coordinates to ``diagonalize'' the map $F$, so that in the new coordinate system, $F = G + E$, where $G$ is an affine map and $E$ is an error term. The map $G$ has the form $(a,b) \mapsto (-2a, b+ c_0)$ and therefore the line $\{ a = 0 \}$ is invariant. Our goal is to show that there is an invariant curve for $F$ as well.
	\item We take a rectangle of the form $X = \{(a,b): R_0 \le b, -R_1 \le a \le R_1 \}$. Assuming that the inverse $F^{-1}$ is well-defined on $X$, we may write $F^{-1} = G^{-1} + H$ and derive conditions $H$ in terms of $R_0$ and $R_1$ which guarantee existence of an invariant curve for $F^{-1}$ (hence for $F$ as well) inside $X$.
	\item We then find a pair $(R_0, R_1)$ for which $F$ is invertible on $X$ and the conditions from (2) on $H$ are satisfied. At this point, we deduce the existence of an invariant curve $\gamma^{inv}$. 
	\item  We show that some forward iterate of $L$ intersects $\gamma^{inv}$ transversally. From this, we obtain the initial value $\alpha_0$. 
	\item  Finally, we show that one can take $R_1 \sim e^{-\mathrm{const} \cdot R_0}$ as $R_0\rightarrow \infty$, which implies that $\gamma^{inv}$ converges exponentially to the vertical line. 
\end{enumerate}

In Subsection \ref{sub:beta=zero}, we give a rather detailed proof in the case $\beta = 0$, following the outline described above. Then in Subsection \ref{sub:beta=positive}, we treat the case $\beta > 0$, but as the structure of the argument is similar, we mainly indicate the necessary modifications. 

\begin{nota}
 	Given a continuous function $f$ on $A$, we write $\V f \V_A = \sup_{x \in A} |f(x)|.$ When $f = (f_{ij})$ is a matrix-valued function, we similarly write $\V f \V_A = \max_{i,j} \sup_{x \in A} |f_{ij}(x)|$. Moreover, given two sequences $c_j$ and $d_j$, we write $c_j \approx d_j$ when $c_j/d_j \rightarrow 1$ as $j \rightarrow +\infty$. 
\end{nota}

\section{Proof of the Main Result}\label{sec:main}

\subsection{The Case of $\beta = 0$}\label{sub:beta=zero}
\begin{step}
Since we are concerned with the region $x, y > 0$, we can make a logarithmic change of variables $ u = \ln x, v = \ln y$. In this coordinate system, the map $F$ takes the form \begin{equation*}
(u,v) \mapsto (v, 2u-v + c_0 + \ln(1+ c_1 e^{-2u+v})),
\end{equation*}
where $c_0 = \ln(\l^2)$ and $c_1 = \l^{-2}$. We then diagonalize the affine part of $F$ by another change of coordinates $a = u- v + c_0/3, b = 2u + v$. In this coordinate system, we have \begin{equation}
F:(a,b) \mapsto (-2a,b+c_0) + (-e(a,b),e(a,b)),
\end{equation}
where $e(a,b) = \ln(1+ c_2 e^{-(4a+b)/3} )$ with $c_2 = \l^{-26/9}$. We define the affine part and the error part of $F$ by \begin{equation*}
\begin{split}
G(a,b) &= (-2a,b+c_0),  \qquad E(a,b) = (-e(a,b),e(a,b)).
\end{split}
\end{equation*}

\end{step} 

\begin{step}
	Take two positive numbers $R_0, R_1$ and consider the rectangles
	\begin{equation*}
	\begin{split}
	X &= [-R_1,R_1] \times [R_0,\infty), \qquad X^+ = [-R_1,R_1] \times [R_0 - R_1 -c_0,\infty)
	\end{split}
	\end{equation*} in the $(a,b)$-plane. Later we will choose $R_0, R_1$ so that $F(X^+) \supset X$, which (together with injectivity of $F$) means that $F^{-1}$ is a well-defined as a map $X \rightarrow X^{+}$. Then, writing $F^{-1} = G^{-1} + H$, we obtain \begin{equation}\label{eq:H_and_E}
	H \circ F = - \nabla G^{-1} \circ E.
	\end{equation}
	We write the components of $H$ and $E$ as $H(a,b) = (h_1(a,b),h_2(a,b))$ and $E(a,b) = (e_1(a,b),e_2(a,b))$.
	
	Consider a space of Lipschitz continuous curves whose images lie in the set $X$ with lipschitz constants not exceeding 1:
	\begin{equation*}
	\Gamma=\{ \gamma:[R_0 ,\infty) \rightarrow [-R_1,R_1] : |\gamma'| \leq 1 \},
	\end{equation*}
	Equipped with the metric \begin{equation*}
	d(\gamma_1,\gamma_2) = \sup_{R_0 \le x} |\gamma_1(x)- \gamma_2(x)|,
	\end{equation*}
	$(\Gamma,d)$ becomes a complete metric space. Assume that we have the following bounds on $H$ 
	\begin{equation}\label{eq:boundH}
	\V H \V_X \le \min\{\frac{R_1}{2},c_0 \},
	\end{equation}
	together with the bounds on $\nabla H$
	\begin{equation}\label{eq:bound_on_grad}
	\V \nabla H \V_X \le \frac{1}{10}.
	\end{equation}
	
	\begin{claim}
		The bounds \eqref{eq:boundH}, \eqref{eq:bound_on_grad} guarantee that $F^{-1}$ induces a contraction mapping in $(\G,d)$.
	\end{claim}
	
	Note that in the case when $H$ is identically zero, it is clear that $F^{-1}$ induces a map on $\Gamma$ by taking the image $F^{-1} ( \gamma)$ and cutting away the piece which does not belong to $X$. Denoting this map by $T$, we have $d(T\gamma_1, T\gamma_2) = d(\gamma_1,\gamma_2)/2$ in this special case. To verify the claim, we will list four conditions which together imply the statement, and then proceed to show that \eqref{eq:boundH}, \eqref{eq:bound_on_grad} imply each condition. We begin with
	
	\begin{cond}
		Given $\gamma \in \Gamma$, the image $F^{-1} \gamma$ is contained in $[-R_1,R_1] \times \mathbb{R}$. 
	\end{cond}
	
	This follows directly from $\V h_1 \V_X \le \frac{R_1}{2}$ of \eqref{eq:boundH}. Next,
	
	\begin{cond}
		The image $F^{-1} \gamma$ is the graph of a function $\beta : [R_0 - r, \infty) \rightarrow [-R_1,R_1]$ for some $r \ge 0$.
	\end{cond}
	
	This time, the condition $h_2 \le c_0 \quad \mbox{on}\quad X$ ensures that the image $F^{-1} \gamma$ has a part belonging to the region $\{a \le R_0 \}$. Then, we only need to exclude the possibility that for $t_1 < t_2$ in $[R_0,\infty)$, $F^{-1}(\gamma(t_1),t_1)$ and $F^{-1}(\gamma(t_2),t_2)$ have the same $a$-coordinate, i.e. $	t_1 + h_2(\gamma(t_1),t_1) = t_2 + h_2(\gamma(t_2),t_2).$ By setting $\Delta t = t_2 - t_1$, we note that above equality implies \begin{equation*}
	\Delta t \le \big( \V \partial_2 h_2 \V_X + \V \partial_1 h_2 \V_X \cdot \vert \gamma' \vert \big) \Delta t.
	\end{equation*}
	This is a contradiction since $\V \partial_2 h_2 \V_X + \V \partial_1 h_2 \V_X \le 1/5$ by \eqref{eq:bound_on_grad}. Therefore, we may cut the piece of $F^{-1} \gamma$ not contained in $X$ and define the resulting curve $[R_0,\infty) \rightarrow [-R_1,R_1]$ as $T\gamma$. To show that the curve obtained in this way belongs to $\G$, we only need to check
	
	\begin{cond}
		The curve $T\gamma$ defined above is continuous with Lipschitz constant not exceeding 1. 
	\end{cond}
	
	We show that \eqref{eq:bound_on_grad} implies \begin{equation}\label{eq:Lip}
	\frac{|-\frac{1}{2} (\gamma(t_2)- \gamma(t_1)) + h_1(\gamma(t_2),t_2) - h_1(\gamma(t_1),t_1) |}{\Delta t + (h_2(\gamma(t_2),t_2) - h_2(\gamma(t_1),t_1)} \le 1,
	\end{equation}
	for any $t_2 > t_1 \ge R_0$ and $\g \in \G$.
	
	Indeed, the denominator is bounded below by \begin{equation*}
	\D t (1 - \V \pr_2 h_2 \V_X - \V \pr_1 h_2 \V_X ) \ge \frac{4}{5} \D t,
	\end{equation*}
	and the first term on the numerator satisfies \begin{equation*}
	\left|-\frac{1}{2} (\gamma(t_2)- \gamma(t_1)) \right| \le \frac{1}{2} \D t,
	\end{equation*}
	while the second one satisfies \begin{equation*}
	|h_1(\gamma(t_2),t_2) - h_1(\gamma(t_1),t_1)| \le (\V \pr_2 h_1 \V_X + \V \pr_1 h_1\V_X) \D t \le \frac{1}{5} \D t.
	\end{equation*}
	
	Combining these, we obtain \eqref{eq:Lip}. That is, $T$ defines a dynamical system on the set $\Gamma$. Finally, we require that
	
	\begin{cond}
		The map $T$ is a contraction on $\G$.
	\end{cond}
	
	We again verify that \eqref{eq:bound_on_grad} is enough to establish it. Take two curves $\gamma_1, \gamma_2 \in \G$ and $t \ge R_0$. Denote $F^{-1}(\gamma_i(t),t) = O_i$, and set $O'$ to be the point on the image $F^{-1}\gamma_2$ which has the same $b$-coordinates with $O_1$. It will be enough to show that $d(O_1, O') \le \mu |\gamma_1(t) - \gamma_2(t)|$ with some $\mu < 1$. For this we will bound each of $d(O_1,O_2)$ and $d(O_2,O')$ in terms of $|\gamma_1(t) - \gamma_2(t)|$.
	
	First, \begin{equation*}
	\begin{split}
	d(O_1,O_2) &\le d(G^{-1}(\gamma_1(t),t) - G^{-1}(\gamma_2(t),t)) +  d(H(\gamma_1(t),t) - H(\gamma_2(t),t)) \\
	&\le \big(\frac{1}{2}+ \sqrt{\V \pr_1h_1 \V_X^2 + \V \pr_1 h_2 \V_X^2} \big) |\gamma_1(t)-\gamma_2(t)|\\
	&\le \big(\frac{1}{2} +\frac{\sqrt{2}}{10} \big) |\gamma_1(t) - \gamma_2(t)|.
	\end{split}	
	\end{equation*}
	
	Next, we set $t^*$ to be the point such that the image of $(\gamma_2(t^*),t^*)$ by $F^{-1}$ has the same $b$-coordinates with $O_1$. Then $t^*$ is determined by the equation \begin{equation}\label{eq:tstar}
	|t^*-t| = |h_2(\gamma_2(t^*),t^*) - h_2(\gamma_1(t),t)|,
	\end{equation}
	
	We bound the right hand side of \eqref{eq:tstar} by 
	\begin{equation*}
	|t^* - t| \le |\pr_1 h_2| \cdot |\gamma_2(t^*) - \gamma_1(t)| + |\pr_2 h_2| \cdot |t^* - t|
	\end{equation*}
	which in turn implies 
	\begin{equation*}
	\begin{split}
	\frac{1}{1- \V \pr_2h_2\V_X } |t^*-t| &\le \V \pr_1h_2 \V_X \cdot |\gamma_2(t^*) - \gamma_1(t)| \\
	&\le \V \pr_1h_2 \V_X \cdot
	\big( |\gamma_2(t^*) - \gamma_2(t)| + |\gamma_2(t) - \gamma_1(t)| \big) \\
	&\le \V \pr_1h_2 \V_X \cdot 
	\big( |t^* -t | + | \gamma_2(t) - \gamma_1(t)| \big),
	\end{split}
	\end{equation*}
	and from \eqref{eq:bound_on_grad}, we obtain \begin{equation*}
	|t^* -t | \le \frac{1}{10} |\gamma_2(t) - \gamma_1(t)|.
	\end{equation*}
	
	Now we can bound $d(O_2,O')^2 \le I^2 + II^2,$ where $I$ and $II$ denote the difference of $O_2,O'$ in $a$ and $b$ coordinates, respectively. We have 
	\begin{equation*}
	\begin{split}
	I &\le \frac{1}{2}|t^* - t| + |h_1 (\gamma_2(t^*),t^*) - h_1(\gamma_2(t),t)| \\
	&\le \big( \frac{1}{2} + \V \pr_1h_1\V_X + \V \pr_2h_1\V_X \big) |t^* - t|,
	\end{split}
	\end{equation*}
	and similarly,
	\begin{equation*}
	\begin{split}
	II &\le |t^* - t| + |h_2 (\gamma_2(t^*),t^*) - h_2(\gamma_2(t),t)| \\
	&\le \big( 1 + \V \pr_1h_2\V_X + \V \pr_2h_2\V_X \big) |t^* - t|,
	\end{split}
	\end{equation*}
	
	With \eqref{eq:bound_on_grad}, we get \begin{equation*}
	d(O_2,O') \le \sqrt{2} |t^* -t| \le \frac{\sqrt{2}}{10} |\gamma_1(t) - \gamma_2(t)|,
	\end{equation*}
	and finally
	\begin{equation*}
	d(O_1,O') \le \big( \frac{1}{2} + \frac{\sqrt{2}}{5} \big) |\gamma_1(t) - \gamma_2(t)|.
	\end{equation*}
	
\end{step}

\begin{step}
	We will pick a pair $(R_0,R_1)$ so that the region $X$ satisfies $F(X^+) \supset X$ and the bounds \eqref{eq:boundH},\eqref{eq:bound_on_grad}. It is easy to see that requiring \begin{equation}\label{eq:bound_on_e}
\V e_1 \V_{X^+}, \V e_2  \V_{X^+} \le \min\{R_1,c_0\}
\end{equation}
guarantees $F(X^+) \supset X$. Next, from the expression $H \circ F (a,b) = - \nabla G^{-1} \circ E = (e_1(a,b)/2,-e_2(a,b))$, we see that \eqref{eq:bound_on_e} is sufficient to guarantee bounds on $H$ in \eqref{eq:boundH}. 

Before we proceed further, let us fix a convention for matrix entries: $( \nabla A )_{ij} = \pr_j A_i$. With this notation, we have \begin{equation}\label{eq:gradH}
(\nabla H \circ F) \nabla F  = \nabla ( -\nabla G^{-1} \circ E ) = \begin{pmatrix}
-\frac{1}{2} \pr_1 e_1 & -\frac{1}{2} \pr_2e_1 \\
\pr_1e_2 & \pr_2e_2
\end{pmatrix}.
\end{equation}

From $\nabla F = \nabla G + \nabla E$, we write \begin{equation*}
\begin{split}
(\nabla F)^{-1} &= \nabla G^{-1} \big( I + (\nabla E \nabla G^{-1} ) \big)^{-1} = \nabla G^{-1} \cdot \big( \sum_{n \ge 0} (\nabla E \nabla G^{-1} )^n \big).
\end{split}
\end{equation*}

Since $\V (\nabla E \nabla G^{-1} ) \V_{X^+} \le \V \nabla E \V_{X^+}$, we have \begin{equation*}
\V (\nabla F)^{-1}\V_{X^+} \le \frac{1}{1-\V \nabla E \V_{X^+} }
\end{equation*}
and we obtain from \eqref{eq:gradH} that \begin{equation}
\V \nabla H \V_{X} \le \frac{2}{1-\V \nabla E \V_{X^+}} \cdot \V \begin{pmatrix}
	-\frac{1}{2} \pr_1 e_1 & -\frac{1}{2} \pr_2e_1 \\
	\pr_1e_2 & \pr_2e_2
	\end{pmatrix} \V_{X^+} \le \frac{2\V \nabla E \V_{X^+}}{1-\V \nabla E \V_{X^+}}.
\end{equation}

Requiring $ \V \nabla E \V_{X^+} \le \frac{1}{25}$ is sufficient to obtain the bound \eqref{eq:bound_on_grad} on $\nabla H$. In conclusion, the following are sufficient conditions on $R_0$ and $R_1$: \begin{equation}\label{eq:restrictions_on_R}
\begin{split}
\V E \V_{X^+} &\le \min\{c_0,R_1\}, \qquad \V \nabla E \V_{X^+} \le \frac{1}{25}.
\end{split}
\end{equation}

We proceed with the explicit formula for the error. Since $|e_i(a,b)| = \ln(1+ c_2 e^{-(4a+b)/3})$ for $ i = 1,2$, we get the maximal value of error in $X^+$ upon substituting $a = R_0 - R_1 - c_0, b = -R_1$. Hence, 
\begin{equation}
\V e_i \V_{X^+} \le c_2 e^{-\frac{1}{3}(4(R_0 - R_1 - c_0) - R_1)} = \lambda^{-\frac{2}{9}} e^{-\frac{1}{3}(4R_0 - 5R_1)}.
\end{equation}
and \begin{equation}
\V \pr_2 e_i \V_{X^+} \le \frac{1}{3}\lambda^{-\frac{2}{9}} e^{-\frac{1}{3}(4R_0 - 5R_1)},
\end{equation}
while $\V \pr_1 e_i \V_{X^+} = 4\V \pr_2 e_i \V_{X^+} $. Therefore, if we pick $R_1 = \min\{\frac{3}{100}, c_0\}$ and then $R_0$ in a way that $\lambda^{-\frac{2}{9}} e^{-\frac{1}{3}(4R_0 - 5R_1)} \le R_1$, all the requirements in \eqref{eq:restrictions_on_R} are satisfied.

Since $\l = 2$, if we fix $R_1 = 3/100$, any value of $R_0$ not less than \begin{equation*}
\frac{3}{4} \big( \ln\frac{100}{3} - \frac{2}{9}\ln 2 + \frac{5}{100} \big) \approx 2.55
\end{equation*}
would work. In particular, we have obtained the existence of a unique $F$-invariant curve $\gamma^{inv}:[2.56,\infty) \rightarrow [-0.03,0.03]$ whose graph lies in $X$. We will soon take $R_0 \rightarrow \infty$, but by an abuse of notation, let us denote the corresponding restriction of the curve by the same letter $\gamma^{inv}$.
\end{step}

\begin{figure}
	\centering
	\includegraphics[scale=0.4]{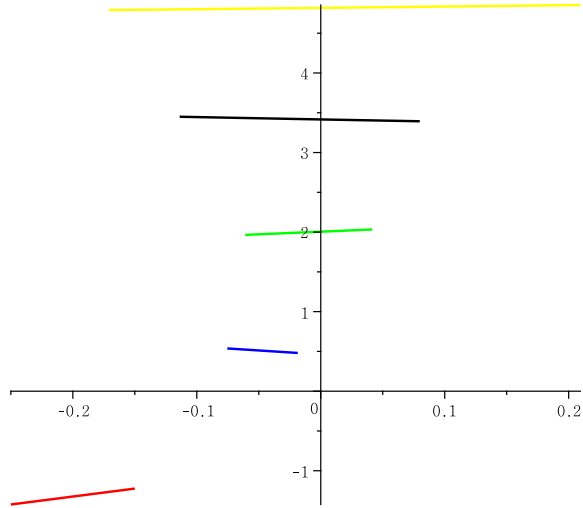}
	\caption{A few forward $F$ iterates of the line segment $I:=\{(t,2t-\frac{2}{3}\ln\l^2): -0.25\le t \le -0.15\}$\label{fig:iterates}: the segments represent $I$, $F(I)$, $F^2(I)$, $F^3(I)$, and $F^4(I)$, from below to above.}
\end{figure}

\begin{step}
	We will show in Lemma \ref{lem:A} that there is an $N>0$ such that $F^N(L)$ intersects $\gamma^{inv}$ (see Figure \ref{fig:iterates}).  By the definition of $L$, we know that the point of intersection has the form $(\alpha_{N-1}, \alpha_N)$ (in the $(x,y)$-coordinates), where $\alpha_{N-1}$ and $\alpha_N$ are from a sequence $\{\alpha_n\}_{n \ge 0}$ satisfying the recurrence \eqref{eq:recursion} with some $\alpha_0 > 0$. We have obtained the value $\alpha_0$. 
\end{step}

\begin{step}
	We now shrink the domain $X$ by taking pairs $(R_0,R_1)$ in a way that $R_0 \rightarrow \infty$ and $R_1 \rightarrow 0$. Note that if we take $R_0$ large then we can take the pair in a way that \begin{equation}
	R_1 \approx e^{-c' R_0}
	\end{equation}
	for some constant $c'>0$. In particular, the curve $\gamma^{inv}$ has the asymptotics \begin{equation}\label{eq:curve_asym}
	|\gamma^{inv}(t)| \le e^{-c't},
	\end{equation}
	as $t \rightarrow \infty$. From the estimate \eqref{eq:curve_asym}, it follows that $\frac{\alpha_n}{\l^{2n/3}} \rightarrow \mathrm{const}$ as $n \rightarrow\infty$. 
\end{step}


\subsection{The Case of $\beta > 0$}\label{sub:beta=positive}

Proceeding analogously as in the case $\beta = 0$, we define the map \begin{equation}\label{eq:F_beta}
F_\beta: (x,y) \mapsto (y, -\frac{\l}{2\beta}y +  \frac{\l}{2\beta}y \sqrt{1 + Z_\beta(x,y)^2}  ),
\end{equation}
where \begin{equation}
Z_\beta(x,y) = \frac{4\beta}{\l}\big( \l^2 \frac{x^2}{y^2} + \beta\l \frac{x}{y} + \frac{1}{y} \big),
\end{equation}
for $0 < \beta < \beta_0$. We may take $\beta_0 = 1$ initially, but we will need to adjust it to be small (and unspecified) at several places from now on.

For convenience, set \begin{equation}
m_\beta(x,y) = -\frac{\l}{2\beta}y +  \frac{\l}{2\beta}y \sqrt{1 + Z_\beta(x,y)^2}.
\end{equation}

It is easy to check that the map $F_\beta$ is well-defined as a map $\mathbb{R}_+ \times \mathbb{R}_+$ to itself and is injective in this region. Indeed, we will consider $F_\beta$ in the region $\{(x,y): x,y>0, Z_\beta(x,y) \le 1/2 \}$ so that we can take the Taylor expansion of $(1+Z_\beta(x,y))^{1/2}$. This will be achieved by restricting to the values of $x,y>0$ with $x/y \le r_0$ and $1/y \le r_1$. Any large values of $r_0, r_1 > 0 $ are allowed at the cost of taking $\beta_0$ small. 

By Taylor expanding $\sqrt{1+Z^2}$ and collecting terms of the same degress in $x$ and $y$, we obtain \begin{equation}\label{eq:expansion_of_beta}
m_\beta(x,y) = \sum_{n \ge 0} d^+_n(\beta) \frac{x^{n+1}}{y^n} + \sum_{n , k \ge 0} d^-_{n,k}(\beta) \frac{x^n}{y^{n+k}},
\end{equation}
which is uniformly convergent for the set of pairs  $(x,y)$ we consider. Notice that the term $-\l y/2\beta$ gets cancelled and all the terms of \eqref{eq:expansion_of_beta} are $O(1)$ as $\beta \rightarrow 0$. With logarithmic change of coordinates $u = \ln x$ and $v = \ln y$, the form of $F_\beta$ becomes $(u,v) \mapsto (v, \tilde{m}_\beta(u,v))$,
with 
\begin{equation}\label{eq:tilde_m}
\begin{split}
\tilde{m}_\beta(u,v) &= 2u - v  +\ln \big( (\sum_{n \ge 0} d^+_{n}(\beta) e^{(n-1)(u-v)}) + (\sum_{n,k\ge 0} d^-_{n,k}(\beta) e^{(n-2)(u-v) -(k+1)v}) \big) \\
&=:2u-v + \ln \big( f_1(\beta)(u,v) + f_2(\beta)(u,v) \big)\\
&= 2u-v + c_0 + \ln \big( 1 +  (f_1(\beta)(u,v)/\l^2 -1) + f_2(\beta)(u,v) \big),
\end{split}
\end{equation}
and one can check that the term $(f_1(\beta)(u,v)/\l^2 -1)$ vanishes on the line $\{u =v\}$. Diagonalizing the affine part by the change of coordinates $
a = u-v + \frac{1}{3}c_0 $ and $b = 2u+v$, we arrive at the form \begin{equation}
F_\beta:(a,b) \mapsto (-2a -e_\beta(a,b), b + c_0 + e_\beta(a,b) ).
\end{equation}
From the expansion in \eqref{eq:tilde_m}, one can check that the error term $e_\beta$ has the form \begin{equation}\label{eq:form_of_error}
e_\beta(a,b) = \ln\big( 1 + g_1(\beta)(a) + g_2(\beta)(a,b)e^{-b/3} \big),
\end{equation}
with estimates \begin{equation}\label{eq:bound_g1}
\begin{split}
|g_1(\beta)(a)| \le \beta C |a|
\end{split}
\end{equation}
and 
\begin{equation}\label{eq:bound_g2}
|g_2(\beta)(a,b) - \l^{-26/9} | \le \beta C 
\end{equation}
in a region of the form $|a| \le \ln r_0, b \ge r_2 = r_2(r_0,r_1)$, for some constant $C = C(r_0,r_2)> 0$, and for $\beta \in [0,\beta_0)$. We can easily deduce bounds of the similar form for the partial derivatives of $e_\beta(a,b)$. 

The point is that, if we consider $X = [-R_1,R_1] \times [R_0,\infty)$ with $R_1 \ll 1$ and $R_0 \gg 1$, then the expansions \eqref{eq:expansion_of_beta}, \eqref{eq:tilde_m} and the bounds \eqref{eq:bound_g1},\eqref{eq:bound_g2} are valid and we have uniform convergence $\V e - e_\beta \V_X \rightarrow 0,$ as $\beta \rightarrow 0$. Therefore, \textbf{Steps 1--3} from the previous section goes through literally in this case as well, with the only difference being that we may need to take $R_1 \ll 1$ and $R_0 \gg 1$. We conclude the existence of the invariant curve $\gamma_\beta^{inv}$, for each $0 < \beta < \beta_0$. By the uniform convergence in $\beta$ of the error term, we deduce that the invariant curves themselves converge uniformly to $\gamma^{inv}$. Since the intersection between $F^N(L)$ and $\gamma^{inv}$ was transversal for all $N$ (whenever they intersect), possibly after taking a smaller value of $\beta_0$, for $0 < \beta < \beta_0$, there is an iterate $F^{N'} L$ of the initial line $L$ which crosses $\gamma^{inv}_\beta$. This takes care of \textbf{Step 4}. At this point, we have obtained the values $\{\alpha_0(\beta)\}_{0<\beta<\beta_0}$. Finally, \textbf{Step 5} follows from the exponential decay of the error in $b \rightarrow \infty$ of \eqref{eq:form_of_error}.
To conclude the proof, it only remains to establish the following 
\begin{lemma}\label{lem:A}
	Let us denote the line segment $\{(t,2t-\frac{2}{3}\ln2^2): -0.4\le t \le 0\}$ by $I$. Then $F^3(I)$ intersects $\gamma^{inv}$ transversally.
\end{lemma}

\begin{proof}
	
	We investigate each iterate of $I$. To begin with, one sees that the image $F(I)$ is indeed a segment of the line $\{ a + b =\ln 2^{\frac{2}{3}}  \}$. An explicit computation, shows that $F(I)$ contains the line segment \begin{equation*}
	J:= \{(s,-s + \ln 2^{\frac{2}{3}}): -0.1 \le s \le -0.01 \}.
	\end{equation*} From now on, we will always assume that the variable $s$ takes values in $[-0.1,-0.01]$.  Parametrizing the set $F(J)$ by $s$, we have \begin{equation*}
	(-2s,-s + \ln 2^{\frac{8}{3}}) + (-e(s),e(s)), \qquad\mbox{with}\qquad e(s) = \ln(1+ 2^{-\frac{10}{3}}e^{-2s}),
	\end{equation*}
 with crude estimates $|e(s)|< 1/8$ and $|e'(s)| < 1/4.$ Next, we parametrize the set $F^2(J)$ as follows: \begin{equation*}
	(4s+2e(s),-s+e(s) + \ln 2^{\frac{14}{3}}) + (-E(s),E(s)),
	\end{equation*} with $E(s) = \ln(1 + 2^{-\frac{42}{9}} e^{6s + 2e(s)}),$
	and we have  $|E(s)| < 1/16$ and $|E'(s)| < 3/8.$ We note that $F^2(J)$ is the graph of a strictly decreasing function defined on $-0.1\le s \le -0.01 $. Indeed, it is enough to check $(4s + 2e(s) - E(s))' > 0$ and $(-s + e(s) + E(s))'<0,$	which follows from the inequalities above. Moreover, from the same estimates, we see that $-s+e(s) + \ln 2^{\frac{14}{3}} + E(s) > 3,$ that is, $F^2(J)$ lies above the line $ b= 3$. Finally, if we plug in the values $ s = -0.1$ and $s = -0.01$, we obtain \begin{equation*}
	4s + 2e(s) - E(s) < -0.03,\qquad  4s + 2e(s) - E(s) > +0.03
	\end{equation*} respectively, using $\ln(1+x) \ge x/2$ for $0<x<1$. This concludes the proof.
\end{proof}

\section*{Acknowledgments} The author thanks Prof. Dong Li for suggesting the problem, and for very helpful discussions. 

\bibliographystyle{plain}
\bibliography{dyadic}

\end{document}